\documentclass{amsart}
\usepackage{latexsym}
\usepackage{amsmath, amssymb}
\usepackage[all]{xy}
\newtheorem{dfs}{Definition}[section]
\newtheorem{lms}[dfs]{Lemma}
\newtheorem{thms}[dfs]{Theorem}
\newtheorem{props}[dfs]{Proposition}

\newtheorem{cors}[dfs]{Corollary}
\newtheorem{rems}[dfs]{Remark}

\newtheorem*{thm*}{Theorem}

\newcommand\VA{\mathrm{V}(A)}
\newcommand\WA{\mathrm{W}(A)}
\newcommand\TA{\mathrm{T}(A)}
\newcommand\TB{\mathrm{T}(B)}
\newcommand\KA{\mathrm{K}_0^*(A)}

\newcommand{\laff}{\ensuremath{\mathrm{LAff}}}
\newcommand{\baff}{\ensuremath{\mathrm{Aff}_b(\TA)}}
\newcommand{\lap}{\mathrm{LAff}_{\mathrm{b}}(\TA)^{++}}
\newcommand\Asa{A_\mathrm{sa}}
\newcommand\p{\varphi}

\newcommand{\Z}{\ensuremath{\mathcal{Z}}}
\newcommand{\dt}{\ensuremath{d_{\tau}}}

\newcommand\lspan{\mathrm{span}\,}

\begin{document}

\title[The Cuntz Semigroup]{The Cuntz Semigroup, the Elliott Conjecture, \\ and dimension functions
on C$^*$-algebras}

\author{Nathanial P. Brown, Francesc Perera and Andrew S. Toms}
\address{Department of Mathematics, Penn State University, State College, PA 16802, USA}
\email{nbrown@math.psu.edu}
\address{Departament de Matem\`atiques, Universitat Aut\`onoma de Barcelona, 08193 Bellaterra, Barcelona, Spain}
\email{perera@mat.uab.cat}
\address{Department of Mathematics and Statistics, York University, 4700 Keele St., Toronto, Ontario, M3J 1P3 Canada}
\email{atoms@mathstat.yorku.ca}
\thanks{N.B.\ was partially supported by DMS-0554870; F.P.\ by the DGI MEC-FEDER through Project MTM2005-00934, and the Comissionat per Universitats i Recerca de la Generalitat de
Catalunya; A.T.\ was supported by an NSERC Discovery Grant.}

\begin{abstract}
We prove that the Cuntz semigroup is recovered functorially from the Elliott invariant for a large
class of C$^*$-algebras.  In particular, our results apply to the largest class of simple
C$^*$-algebras for which K-theoretic classification can be hoped for.  This work has three
significant consequences.  First, it provides new conceptual insight into Elliott's classification
program, proving that the usual form of the Elliott conjecture is equivalent, among $\Z$-stable
algebras, to a conjecture which is in general substantially weaker and for which there are no
known counterexamples.  Second and third, it resolves, for the class of algebras above, two
conjectures of Blackadar and Handelman concerning the basic structure of dimension functions on
C$^*$-algebras.  We also prove in passing that the Cuntz-Pedersen semigroup is recovered
functorially from the Elliott invariant for a large class of simple unital C$^*$-algebras.
\end{abstract}

\maketitle

\section{Introduction}

The Cuntz semigroup $\WA$ of a C$^*$-algebra $A$ is an analogue for positive elements of the semigroup
of Murray-von Neumann equivalence classes of projections $\VA$.  It is deeply connected to the classification
program for simple separable nuclear C$^*$-algebras:  such algebras cannot be classified up to isomorphism
by their $\mathrm{K}$-theory and traces if the natural partial order on the Cuntz semigroup is not determined
by traces in a weak sense (\cite{T2}), and the converse --- known to hold in some cases (\cite{TW3}) ---
may well hold in great generality.  One thus expects the structure of the Cuntz semigroup to be implicit
in the $\mathrm{K}$-theory and traces of sufficiently well-behaved C$^*$-algebras, despite the fact that
computing the Cuntz semigroup of an Abelian C$^*$-algebra is totally infeasible (see \cite[Lemma 5.1]{T2}).

The largest class of simple separable nuclear C$^*$-algebras which one may hope will be classified
by the Elliott invariant of $\mathrm{K}$-theoretic data consists of those algebras which absorb
the Jiang-Su algebra $\mathcal{Z}$ tensorially.  Such algebras are said to be {\it
$\mathcal{Z}$-stable}.  It is expected, but not yet known, that simple unital approximately
homogeneous (AH) algebras of slow dimension growth will be $\mathcal{Z}$-stable.  Before stating
the first of our two main results, we recall that the Elliott invariant of a unital C$^*$-algebra
$A$ is the 4-tuple
\[
\mathrm{Ell}(A) := \left( (\mathrm{K}_0A,\mathrm{K}_0A^+,[1_A]),\mathrm{K}_1A,\mathrm{T}A,\rho_A
\right),
\]
where the $\mathrm{K}$-groups are the Banach algebra ones, $\mathrm{T}A$ is the tracial state
space, and $\rho_A$ is the pairing between $\mathrm{K_0}$ and $\mathrm{T}A$ given by evaluating a
$\mathrm{K}_0$-class at a trace.

\vspace{2mm}
\noindent
{\bf Theorem A.} {\it Let $A$ be a simple unital finite C$^*$-algebra which is either
exact and $\mathcal{Z}$-stable or AH of slow dimension growth.  Then, there is a functor
which recovers $\WA$ from the Elliott invariant $\mathrm{Ell}(A)$.}

\vspace{2mm} The functor of Theorem A is defined in Section 2, and describes $\WA$ in terms of the
Murray-von Neumann semigroup $\VA$ and certain affine functions on the tracial state space
$\mathrm{T}(A)$ of $A$.  It computes $\WA$ for the majority of our stock-in-trade simple separable
nuclear C$^*$-algebras.

The usual form of Elliott's classification conjecture states that isomorphisms between
the Elliott invariants of simple separable nuclear C$^*$-algebras are liftable to
isomorphisms between the algebras.  An immediate consequence of Theorem A is that Elliott's
conjecture is equivalent, among $\mathcal{Z}$-stable algebras, to a formally weaker conjecture:
isomorphisms at the level of the invariant $(\mathrm{Ell}(\bullet),\mathrm{W}(\bullet))$ are
liftable to isomorphisms at the level of algebras.  Outside the class of $\mathcal{Z}$-stable
algebras, there are no known counterexamples to this weaker conjecture.  Thus, we reconcile
the necessity of the Cuntz semigroup in any effort to classify general separable nuclear
C$^*$-algebras with its absence from known classification results, and provide
a new point of departure for proving classification theorems for $\mathcal{Z}$-stable algebras.
It should be noted that the Cuntz semigroup stands out among candidates for augmenting the
Elliott invariant.  Even the addition of all continuous homotopy invariant functors and all
stable invariants to $\mathrm{Ell}(A)$ does not yield a complete invariant (\cite{T3}).

Our second main result concerns two conjectures of Blackadar and Handelman on the structure
of dimension functions on C$^*$-algebras, put forth in the early 1980s (see Section \ref{sec:prelim} for
terminology):
\begin{enumerate}
\item[(i)] the lower semicontinuous dimension functions  are dense in the space of all dimension functions;
\item[(ii)] the affine space of dimension functions is a simplex.
\end{enumerate}
We will refer to these as the Blackadar-Handelman conjectures for brevity.
Blackadar and Handelman proved that (i) holds for commutative $C^*$-algebras, but did not prove
(ii) in any case. The only further progress on these conjectures was made by the second
named author in~\cite[Corollary 4.4]{P1}, where (ii) was confirmed for the class
of unital $C^*$-algebras with real rank zero and stable rank one.  We obtain:

\vspace{2mm}
\noindent
{\bf Theorem B.}  {\it Let $A$ be a simple unital finite C$^*$-algebra which is either
exact and $\mathcal{Z}$-stable or AH of slow dimension growth.  Then, the Blackadar-Handelman
conjectures hold for $A$, and the simplex of dimension functions on $A$ is in fact a Choquet
simplex.}

\vspace{2mm}
Theorem B applies, in particular,
to several classes of $\Z$-stable ASH algebras (see Section 4 of \cite{TW2}).  The $\Z$-stability of
these algebras can only be established using the fact that they satisfy the Elliott conjecture.
Thus, our confirmation of the Blackadar-Handelman conjectures for these algebras constitutes
a bona fide application of $\mathrm{K}$-theoretic classification theorems.

The proofs of Theorems A and B use the fact that the C$^*$-algebras in question have strict
comparison of positive elements in a crucial way.  (This property says roughly that the order
structure on the Cuntz semigroup is determined by traces---see Section 2 for an exact definition.)
For the $\mathcal{Z}$-stable algebras we consider, this property was established by M. R{\o}rdam
in \cite{R1}, building on his earlier work in \cite{Rfunct}.  For our AH algebras, strict
comparison was established by the third named author in \cite{tomscomp}.

\vspace{2mm}

{\noindent\bf Acknowledgement:} Most of this work was carried out while the first and second named
authors visited the third at the University of New Brunswick.  We gratefully acknowledge the
support and hospitality extended by UNB during our visit.

%%%%%%%%%%%%%%%%%%%%%%%%%%%%%%%%%%%%%%%%%%%%%
\section{Preliminaries and notation}
\label{sec:prelim}
\subsection*{Cuntz Equivalence}

From here on we make the blanket assumption that all
C$^*$-algebras are separable unless otherwise stated or obviously false.

Let $A$ be a C$^*$-algebra, and let $\mathrm{M}_n(A)$ denote the $n \times n$
matrices whose entries are elements of $A$.  If $A = \mathbb{C}$, then we simply write $\mathrm{M}_n$.
Let $\mathrm{M}_{\infty}(A)$ denote the algebraic limit of the
direct system $(\mathrm{M}_n(A),\phi_n)$, where $\phi_n:\mathrm{M}_n(A) \to \mathrm{M}_{n+1}(A)$
is given by
\[
a \mapsto \left( \begin{array}{cc} a & 0 \\ 0 & 0 \end{array} \right).
\]
Let $\mathrm{M}_{\infty}(A)_+$ (resp. $\mathrm{M}_n(A)_+$)
denote the positive elements in $\mathrm{M}_{\infty}(A)$ (resp. $\mathrm{M}_n(A)$).
For positive elements $a$ and $b$ in $\mathrm{M}_{\infty}(A)$, write $a\oplus b$ to denote the element $\left( \begin{array}{cc} a & 0 \\ 0 & b \end{array} \right)$, which is also positive in $\mathrm{M}_{\infty}(A)$.

Given $a,b \in \mathrm{M}_{\infty}(A)_+$, we say that $a$ is {\it Cuntz subequivalent} to
$b$ (written $a \precsim b$) if there is a sequence $(v_n)_{n=1}^{\infty}$ of
elements of $\mathrm{M}_{\infty}(A)$ such that
\[
\Vert v_nbv_n^*-a\Vert \stackrel{n \to \infty}{\longrightarrow} 0.
\]
We say that $a$ and $b$ are {\it Cuntz equivalent} (written $a \sim b$) if
$a \precsim b$ and $b \precsim a$.  This relation is an equivalence relation,
and we write $\langle a \rangle$ for the equivalence class of $a$.  The set
\[
\WA := \mathrm{M}_{\infty}(A)_+/ \sim
\]
becomes a positively ordered Abelian monoid when equipped with the operation
\[
\langle a \rangle + \langle b \rangle = \langle a \oplus b \rangle
\]
and the partial order
\[
\langle a \rangle \leq \langle b \rangle \Leftrightarrow a \precsim b.
\]
In the sequel, we refer to this object as the {\it Cuntz semigroup} of $A$.

Given $a$ in $\mathrm{M}_{\infty}(A)_+$ and $\epsilon > 0$, we denote by
$(a-\epsilon)_+$ the element of $C^*(a)$ corresponding (via the functional
calculus) to the function
\[
f(t) = \mathrm{max}\{0,t-\epsilon\}, \ t \in \sigma(a).
\]
(Here $\sigma(a)$ denotes the spectrum of $a$.)

In order to ease the notation, we will let
$A_a$ denote the hereditary C$^*$-algebra generated by a
positive element $a$ in $A$, that is, $A_a=\overline{aAa}$. Recall
that if $A$ is a separable C$^*$-algebra, then all hereditary
algebras are of this form.  Some of our results require the assumption
that $A$ has stable rank one, that is, the set of invertible elements
is dense in $A$.  We denote the stable rank of $A$ by $\mathrm{sr}(A)$.
If $\mathrm{sr}(A)=1$, then Cuntz subequivalence, viewed as a relation
on hereditary subalgebras, is implemented by unitaries (\cite{Rfunct}).

The proposition below collects some facts about Cuntz subequivalence due
to Kirchberg and R{\o}rdam.

\begin{props}[Kirchberg-R{\o}rdam (\cite{KR}), R{\o}rdam (\cite{Rfunct})]\label{basics}
Let $A$ be a C$^*$-algebra, and $a,b \in A_+$.
\begin{enumerate}
\item[(i)] $(a-\epsilon)_+ \precsim a$ for every $\epsilon > 0$.
\item[(ii)] The following are equivalent:
\begin{enumerate}
\item[(a)] $a \precsim b$;
\item[(b)] for all $\epsilon > 0$, $(a-\epsilon)_+ \precsim b$;
\item[(c)] for all $\epsilon > 0$, there exists $\delta > 0$ such that $(a-\epsilon)_+ \precsim (b-\delta)_+$.
\end{enumerate}
\item[(iii)] If $\epsilon>0$ and $\Vert a-b\Vert<\epsilon$, then $(a-\epsilon)_+ \precsim b$.
\item[(iv)] If moreover $\mathrm{sr}(A)=1$, then
\[
a\precsim b \text{ if and only if for every }\epsilon>0\,,\text{
there is }u\text{ in }U(A)\text{ such that }u^*(a-\epsilon)_+u\in
A_b\,.
\]
\end{enumerate}
\end{props}

Note that, if $A_a\subseteq A_b$ for positive elements $a$ and $b$ in $A$, we have that $a\precsim b$
(by Proposition~\ref{basics}).

\subsection*{Traces, Quasitraces, States and Dimension Functions}

As usual, we shall denote the state space of $A$ (that is, the space of positive, unital, linear functionals) by $\mathrm{S}(A)$. The set of tracial states will be denoted by $\TA$ and $\mathrm{QT}(A)$ will be used for
the space of normalised 2-quasitraces on $A$ (v. \cite[Definition II.1.1]{bh}). Note that $\TA\subseteq\mathrm{QT}(A)$, and equality holds when $A$ is exact (see~\cite{Ha}).

Let $\mathrm{St}(\WA,\langle 1_A\rangle)$ denote the set of additive and order preserving maps $s$
from $\WA$ to $\mathbb{R}^+$ having the property that $s(\langle 1_A \rangle) = 1$. Such maps are
generally called \emph{states} and in the particular case of a C$^*$-algebra, they are termed
\emph{dimension functions}.  The set of all dimension functions on $A$ is denoted by
$\mathrm{DF}(A)$.

Given $\tau$ in $\mathrm{QT}(A)$, one may
define a map $d_{\tau}\colon\mathrm{M}_{\infty}(A)_+ \to \mathbb{R}^+$ by
\begin{equation}\label{ldf}
d_{\tau}(a) = \lim_{n \to \infty} \tau(a^{1/n}).
\end{equation}
This map is lower semicontinous, and depends only on the Cuntz equivalence class
of $a$.  It moreover has the following properties:
\vspace{2mm}
\begin{enumerate}
\item[(i)] if $a \precsim b$, then $d_{\tau}(a) \leq d_{\tau}(b)$;
\item[(ii)] if $a$ and $b$ are mutually orthogonal, then $d_{\tau}(a+b) = d_{\tau}(a)+d_{\tau}(b)$;
\item[(iii)] $d_{\tau}((a-\epsilon)_+) \nearrow d_{\tau}(a)$ as $\epsilon \to 0$.
\end{enumerate}
\vspace{2mm}
Thus, $d_{\tau}$ defines a state on $\WA$.
Such states are called \emph{lower semicontinuous dimension functions}, and the set of them
is denoted $\mathrm{LDF}(A)$. It was proved in~\cite[Theorem II.4.4]{bh} that $\mathrm{QT}(A)$ is a simplex;
the map from $\mathrm{QT}(A)$ to $\mathrm{LDF}(A)$ defined by (\ref{ldf}) is
bijective and affine (\cite[Theorem II.2.2]{bh}), but generally not continuous.
We also have that $\mathrm{LDF}(A)$ is a (generally proper) face of $\mathrm{DF}(A)$, see~\cite[Proposition II.4.6]{bh}. If $A$ has the property that $a \precsim b$ whenever $s(a) < s(b)$ for every $s \in \mathrm{LDF}(A)$, then we say that $A$ has
{\it strict comparison of positive elements} or simply {\it strict comparison}.

The Grothendieck group of $\WA$ is denoted by $\KA$. The class of an element $a$ from $\mathrm{M}_{\infty}(A)_+$ will be denoted by $[a]$. This is a partially ordered Abelian group with positive cone $\KA^+=\{[a]-[b]\mid b\precsim a\}$. Observe then that $\mathrm{DF}(A)=\mathrm{St}(\KA,\KA^+,[1_A])$, which is the set of group morphisms $s\colon\KA\to\mathbb{R}$ such that $s(\KA^+)\subseteq\mathbb{R}^+$ and $s([1_A])=1$.

\subsection*{The reconstruction functor}
\label{reconstruct}

Let $A$ be a unital and stably finite C$^*$-algebra
with tracial state space $\mathrm{T}(A)$, and let $\mathrm{LAff}_b(\mathrm{T}(A))^{++}$
denote the set of bounded, strictly positive, lower semicontinuous, and affine
functions on $\mathrm{T}(A)$.   Define a semigroup structure on the disjoint union
\[
\widetilde{W}(A) := V(A) \sqcup \mathrm{LAff}_b(\mathrm{T}(A))^{++}
\]
as follows:
\vspace{2mm}
\begin{enumerate}
\item[(i)] if $x,y \in V(A)$, then their sum is the usual sum in $V(A)$;
\item[(ii)] if $x,y \in \mathrm{LAff}_b(\mathrm{T}(A))^{++}$, then their sum
is the usual (pointwise) sum in $\mathrm{LAff}_b(\mathrm{T}(A))^{++}$;
\item[(iii)] if $x \in V(A)$ and $y \in \mathrm{LAff}_b(\mathrm{T}(A))^{++}$,
then their sum is the usual (pointwise) sum of $\hat{x}$ and $y$ in
$\mathrm{LAff}_b(\mathrm{T}(A))^{++}$, where $\hat{x}(\tau) = \tau(x)$,
$\forall \tau \in \mathrm{T}(A)$.
\end{enumerate}
\vspace{2mm}
Equip $\widetilde{W}(A)$ with the partial order $\leq$ which restricts to the
usual partial order on each of $V(A)$ and $\mathrm{LAff}_b(\mathrm{T}(A))^{++}$,
and which satisfies the following conditions for $x \in V(A)$ and $y \in
\mathrm{LAff}_b(\mathrm{T}(A))^{++}$:
\vspace{2mm}
\begin{enumerate}
\item[(i)] $x \leq y$ if and only if $\hat{x}(\tau) < y(\tau)$, $\forall \tau \in \mathrm{T}(A)$;
\item[(ii)] $y \leq x$ if and only if $y(\tau) \leq \hat{x}(\tau)$, $\forall \tau \in \mathrm{T}(A)$.
\end{enumerate}
\vspace{2mm}
The map $A \mapsto \widetilde{W}(A)$ is a functor --- see Section 4 of \cite{pt} for details.

Now let $A$ be simple.  There is a canonical map
$\phi:\WA \to \widetilde{\mathrm{W}}(A)$ first defined in \cite{pt}.
Let us recall its definition.
Denote by $A_{++}$ those elements of $A_+$ which are not Cuntz equivalent to
a projection in $\mathrm{M}_{\infty}(A)$, and set
\[
W(A)_+ = \{ \langle a \rangle \in W(A) \ | \ a \in \mathrm{M}_{\infty}(A)_{++} \}.
\]
The elements of $A_{++}$ are called {\it purely positive}.  If $A$ has stable rank one, then
$W(A)$ is the disjoint union of $V(A)$ (identified with its image in $W(A)$ via the natural map
$[p]\mapsto\langle p\rangle$), and $W(A)_+$. As observed in~\cite[Corollary 2.9]{pt}, if $A$ is
either simple and stably finite or of stable rank one, we have that $W(A)_+$ is actually a
subsemigroup of $W(A)$, and is absorbing in the sense that $a+b\in W(A)_+$ whenever $a\in V(A)$
and $b\in W(A)_+$.

Let $\iota\colon W(A)_+ \to \mathrm{LAff}_b(\mathrm{T}(A))^{++}$ be
given by $\iota(x)(\tau) := d_{\tau}(x)$.
If $A$ is exact and has strict comparison, then $\iota$ is an order embedding on $W(A)_+$
(\cite[Proposition 3.3]{pt}).
Let $\phi:\WA \to \widetilde{\mathrm{W}}(A)$ be given by $\phi|_{V(A)} =
\mathrm{id}_{V(A)}$ and $\phi|_{W(A)_+} = \iota$.  It is proved in \cite{pt}
that $\phi$ is both everywhere-defined and well-defined.

\begin{thms}[P-T, \cite{pt}, Theorem 4.4]\label{embed}
Let $A$ be a simple, unital, exact, and stably finite $C^*$-algebra
with strict comparison of positive elements.  Then,
\[
\phi:W(A) \to \widetilde{W}(A)
\]
is an order embedding.
\end{thms}
\noindent

Thus, under the hypotheses of Theorem \ref{embed}, $\phi$ is an isomorphism
whenever $\iota$ is surjective.

%%%%%%%%%%%%%%%%%%%%%%%%%%%%%%%%%%%%%%%%

\section{Duality and Traces}
\label{sec:duality}

If $\mathrm{S}(A)$ is the state space of a unital C$^*$-algebra $A$ and
$X=\mathrm{span}_\mathbb{R}\mathrm{S}(A)$ is the $\mathbb{R}$-Banach space of self-adjoint
functionals on $A$ then we have two natural dualities: $$X = (A_{\mathrm{sa}})^* \ \mathrm{and} \
X^* = A^{**}_\mathrm{sa},$$ where $A_{\mathrm{sa}}$ (resp.\  $A^{**}_\mathrm{sa}$) denotes the
self-adjoint elements in $A$ (resp.\ in the enveloping von Neumann algebra $A^{**}$). Kadison's
function representation (cf.\ \cite[Theorem 3.10.3]{pedersen}) is a well-known application of these
two facts:  If $f\colon \mathrm{S}(A) \to \mathbb{R}$ is a bounded affine function then there
exists a unique $T \in A^{**}_\mathrm{sa}$ such that $f(\varphi) = \varphi(T)$, for all $\varphi
\in  \mathrm{S}(A)$, and $\|T\| = \|f\|$; moreover, $f$ is continuous if and only if $T \in \Asa$.

The purpose of this section is to prove analogous results when $\mathrm{S}(A)$ is replaced by $\TA$ (cf.\ Theorem \ref{thm:duality} and Corollary \ref{cor:realizingaffinefunctions}).

\subsection*{Some Conventions and Basics}

For a unital C$^*$-algebra $A$, \emph{we always consider $\TA$ a compact topological space, endowed with the weak-$*$ topology coming from $A^*$.} (Hence, a ``continuous" function on $\TA$ means continuous with respect to this topology.)  We regard the $\mathbb{R}$-linear space $\lspan_\mathbb{R} \TA$ as a $\mathbb{R}$-Banach space, equipped with the restriction of the norm on $A^*$; when thinking of $\lspan_\mathbb{R} \TA$ as a locally convex space with respect to the weak-$*$ topology, we will make this point explicit.  The following proposition is well known.

\begin{props}[Jordan Decomposition]
\label{prop:statesandtraces}
For any unital C$^*$-algebra $A$ and self-adjoint functional $\p \in A^*$, there exist (unique) orthogonal central projections $P_+$, $P_- \in A^{**}$ such that $\p_+(a) :=  \p(aP_+)$ and $\p_-(a) := -\p(aP_-)$ are positive linear functionals, $\p = \p_+ - \p_-$ and $\|\p\| = \|\p_+\| + \|\p_-\|.$

If $\p$ has the property that $\p(a^* a) = \p(aa^*)$ for all $a \in A$ and $\p = \p_+ - \p_-$ is its Jordan decomposition then $\p$, $\p_+$ and $\p_-$ are all tracial functionals. Consequently,
\begin{align*}
\lspan_\mathbb{R} \TA & = \{\varphi \in A^*:\p(a^*) = \p(a)^* \ \mathrm{and} \ \p(a^*a)
= \p(aa^*), \forall a \in A\}\\
& =  \{t_1\tau_1 - t_2\tau_2 : t_i \geq 0, \tau_i \in \TA, i=1,2 \}.
\end{align*}
\end{props}

The proof of the next proposition is well known to anyone familiar with the proof
of Kadison's function representation (cf. \cite{pedersen}).

\begin{props}
\label{prop:affineextension} Let $f\colon \TA \to V$ be an affine function into a
$\mathbb{R}$-vector space $V$.  Then, $f$ has a unique extension to a linear function
$\tilde{f}\colon \lspan_\mathbb{R} \TA \to V$.  If $V$ is a topological vector space and $f$ is
continuous then $\tilde{f}$ is also continuous (with respect to the weak-$*$ topology).
\end{props}

\subsection*{Cuntz-Pedersen Equivalence}

There is another notion of equivalence that one can consider in $A_+$, first studied by Cuntz and
Pedersen in \cite{CP}.  Namely, for positive elements $a,b \in A_+$, we write  $a \sim_{CP} b$ if
there exist elements $u_n \in A$ such that $$a = \sum_{n} u_n^*u_n \ \ \mathrm{and} \ \ b = \sum_n
u_n u_n^*,$$ where convergence is in norm. By~\cite[Proposition 1.1]{pedersen:III}
(see also~\cite[Corollary 3.6]{pedersen:IV}), $\sim_{CP}$ is an
equivalence relation. It follows from the definition, and a change of index set, that if
$a_1 \sim_{CP} b_1$ and $a_2 \sim_{CP} b_2$ then $a_1 + a_2 \sim_{CP} b_1 + b_2$.
Thus we can define the \emph{Cuntz-Pedersen semigroup} to be $A_+$ modulo the equivalence relation $\sim_{CP}$. More generally, $$A_0 = \{a -
b: a,b \in A_+, a\sim_{CP} b\}$$ is a $\mathbb{R}$-linear subspace of $\Asa$.  In fact, \cite[Theorem
2.6]{CP} asserts that $A_0$ is a \emph{norm-closed} subspace, and hence we can factor it out.

\begin{dfs}
\label{defn:Aq} Define a $\mathbb{R}$-Banach space by
$$A^q = A_\mathrm{sa}/A_0,$$ and let $A^q_+ =
q(A_+)$ be the image of $A_+$ under the quotient map $q\colon A_\mathrm{sa} \to A^q$.
\end{dfs}

Since $\lspan_\mathbb{R} \TA$ is a weak-$*$ closed subspace of $(\Asa)^*$, it has a predual
(namely, the quotient of $\Asa$ by the pre-annihilator). More precisely:

\begin{props}\cite[Proposition 2.7]{CP}
\label{prop:Aqpredual}
The dual space of $A^q$ is isometrically isomorphic to $\lspan_\mathbb{R} \TA$. Moreover, the induced weak-$*$ topology agrees with the canonical weak-$*$ topology (coming from $A^*$).
\end{props}

For simple, unital C$^*$-algebras the following description of $A^q_+\setminus \{0\}$ is useful.

\begin{props} \cite[Theorem 5.2 and Corollary 6.4]{CP}
\label{prop:cone}
If $A$ is simple, unital and has at least one tracial state, then $$A^q_+\setminus \{0\} = \{x \in A^q: \tau(x) > 0, \forall \tau \in \TA\}.$$ In addition, $A^q_+$ is isomorphic, as an additive semigroup, to the Cuntz-Pedersen semigroup $A_+/\sim_{CP}$.
\end{props}

\begin{cors}  If unital C$^*$-algebras $A$ and $B$ have non-empty affinely homeomorphic tracial state spaces then $A^q \cong B^q$.

If $A$ and $B$ are simple then the Cuntz-Pedersen semigroups $A_+/\sim_{CP}$ and  $B_+/\sim_{CP}$ are also isomorphic.
\end{cors}

\begin{proof}  Assume that $\TA$ and  $\TB$ are affinely homeomorphic (with respect to the restrictions of the weak-$*$ topologies on $A^*$ and, respectively, $B^*$). Then, thanks to Propositions \ref{prop:affineextension} and \ref{prop:Aqpredual}, $\lspan_\mathbb{R} \TA$ is isomorphic to $\lspan_\mathbb{R} \TB$ as locally convex spaces with respect to the weak-$*$ topologies coming from $A^q$ and, respectively, $B^q$.  Thus their dual spaces --- i.e.\ $A^q$ and $B^q$ --- must be isomorphic too.

It is clear that the induced isomorphism $A^q \to B^q$ will map the set  $\{x \in A^q: \tau(x) > 0, \forall \tau \in \TA\}$ bijectively onto $\{y \in B^q: \tau(y) > 0, \forall \tau \in \TB\}$. It follows that $A^q_+$ will get mapped bijectively onto $B^q_+$; hence Proposition \ref{prop:cone} implies that $A_+/\sim_{CP}$ and  $B_+/\sim_{CP}$ are isomorphic too.
\end{proof}

\subsection*{Tracial Analogue of Kadison's Function Representation}

With the canonical predual of $\lspan_\mathbb{R} \TA$ in hand, our tracial version of Kadison's function representation is within sight.  We just need the dual space.  This has a simple description in terms of the enveloping von Neumann algebra $A^{**}$ (indeed, it may be known to some experts, but we are unaware of a reference).

\begin{lms}\label{lem:Zsa-spanTA} Let $A$ be a unital C$^*$-algebra with tracial simplex $\TA$.
If $Z$ denotes the center of the maximal finite summand of $A^{**}$ then there is an isometric
identification $$Z_\mathrm{sa} = (\lspan_\mathbb{R} \TA)^*.$$ In other words, the predual of $Z$ is
equal to $\lspan_\mathbb{C}\TA.$
\end{lms}

\begin{proof} For convenience, let $M$ denote the maximal finite
summand of $A^{**}$ (hence $A^{**} = M \oplus N$, with $N$ infinite). Let $\Phi\colon M \to Z$ be
the canonical center-valued trace (cf.\ \cite[Theorem 2.4.6]{sakai}); that is, $\Phi$ is a
$\sigma$-weakly continuous faithful conditional expectation onto $Z$ with the property that
$\Phi(xy) = \Phi(yx)$, for all $x,y \in M$, and $\tau \circ \Phi = \tau$ for every tracial state on
$M$. Though a slight abuse of notation, we will let $\Phi(a) \in Z$, $a \in A$, denote the
composition of the maps $$A \hookrightarrow A^{**} \to M \stackrel{\Phi}{\to} Z,$$ where $A^{**}
\to M$ is the canonical quotient map.

For each $\tau \in \TA$ we use the same symbol to denote the normal extension to $A^{**}$. Note
that each such $\tau$ is supported on $M$ --- i.e.\ $\tau|_N = 0$, by maximality of $M$ --- and thus
we have a natural inclusion $\lspan_\mathbb{C}\TA \subset M_*$. Since $Z$ is a subalgebra of $M$,
we have a (linear) restriction map $\lspan_\mathbb{C}\TA \to Z_*$. It is evidently isometric (hence
injective) since
\begin{align*}
\|\tau|_{Z}\|_{Z_*} \leq \|\tau\|_{(A^{**})_*} & = \|\tau\|_{A^*}\\
& = \sup \{|\tau(a)|: a \in
A, \|a\|\leq 1\}\\
& = \sup \{|\tau(\Phi(a))|: a \in A, \|a\|\leq 1\}\\
&\leq \|\tau|_{Z}\|_{Z_*}.
\end{align*}

To prove surjectivity of the restriction map $\lspan_\mathbb{C}\TA \to Z_*$, it suffices to show
that every normal state on $Z$ is the restriction of some tracial state.  So, fix a normal state
$\varphi \in Z_*$ and define a trace $\tau$ on $A$ by $\tau(a) = \varphi(\Phi(a))$. (This is
tracial since $\Phi$ is.) One easily checks that (the normal extension of) $\tau$ restricts to
$\varphi$, using the fact that $\Phi$ is a $\sigma$-weakly continuous conditional expectation. This
establishes the canonical isometric identification $Z_* \cong \lspan_\mathbb{C}\TA$.

It follows that $Z_\mathrm{sa} = (\lspan_\mathbb{R} \TA)^*$, because $Z_\mathrm{sa}$ can be identified with the dual of the self-adjoint, normal functionals on $Z$ --- i.e.\ the dual of $\lspan_\mathbb{R} \TA$.
\end{proof}

Summarizing our duality results, we have:

\begin{thms}
\label{thm:duality} Let $A$ be a unital C$^*$-algebra with tracial simplex $\TA$.
Then  $$\mathrm{span}_\mathbb{R}\TA = (A^q)^* \ \mathrm{and} \ (\mathrm{span}_\mathbb{R}\TA)^* =
Z_\mathrm{sa},$$ where $Z$ denotes the center of the maximal finite summand of $A^{**}$.
\end{thms}

\begin{dfs} For a unital C$^*$-algebra $A$, let $\baff$ denote
the set of $\mathbb{R}$-valued bounded affine functions on $\TA$.
Let $A^{**}_\mathrm{sa} \to \baff$ be the restriction of Kadison's function representation
to the tracial state space: $a\mapsto \hat{a}$, where
$$\hat{a}(\tau) = \tau(a),$$ for all $a \in A^{**}_\mathrm{sa}$ and $\tau \in \TA$.
\end{dfs}

We assume below that $A$ is a unital C$^*$-algebra with at least one tracial state.

\begin{cors}
\label{cor:realizingaffinefunctions}
The mapping $a\mapsto \hat{a}$ gives a linear, order preserving, isometric
identification of $Z_\mathrm{sa}$ with $\baff$.  Moreover, for every continuous
$f \in \baff$ there exists $a \in \Asa$ such that $f(\tau) = \tau(a)$, for all
$\tau \in \TA$; if $A$ is simple and $f(\tau) > 0$, for all $\tau \in \TA$, then
we can find a positive $a \in A_+$ with $f(\tau) = \tau(a)$.
\end{cors}

\begin{proof}  Since $\TA$ is identified with the normal states on $Z$, the mapping $a \mapsto \hat{a} \in
\baff$ is easily seen to be an order preserving, isometric injection of $Z_\mathrm{sa}$ into $\baff$.
(Or, it follows from Kadison's function representation, applied to $Z$, and the fact that $\TA$ is dense
in the set of all states on $Z$.)  Surjectivity follows easily from Proposition \ref{prop:affineextension}
and Lemma \ref{lem:Zsa-spanTA}.

Similarly, if $f \in \baff$ is continuous, then Proposition \ref{prop:affineextension} says we can extend
it to a weak-$*$ continuous linear functional $\tilde{f}$ on $\mathrm{span}_\mathbb{R}\TA$.  Since the
dual of $\mathrm{span}_\mathbb{R}\TA$ with respect to this topology is $A^q$, and $A^q$ is a quotient
of $\Asa$, we simply identify $\tilde{f}$ with an element in $A^q$ and lift it to $\Asa$.

When $A$ is simple, every $x\in A^q$ with the property that $\tau(x) > 0$, for all $\tau \in \TA$, can
be lifted to a positive element thanks to Proposition \ref{prop:cone}. This implies the last statement,
so the proof is complete.
\end{proof}

%%%%%%%%%%%%%%%%%%%%%%%%%%%%%%%%%%%%%%%%
\section{Suprema in the Cuntz semigroup}
\label{sec:suprema}

In this section we prove that for C$^*$-algebras with stable rank
one, the Cuntz semigroup admits suprema of countable bounded
sequences in a sense that we now proceed to define.

\begin{dfs}
Let $(M,\leq)$ be a preordered Abelian semigroup with identity
element $0$. We say that an element $x$ in $M$ is the
\emph{supremum} of an increasing sequence $(x_n)$ of elements in
$M$ provided that $x_n\leq x$ for each $n$ and is the smallest
such $x$, meaning that if $y\in M$ and $x_n\leq y$ for all $n$,
then necessarily $x\leq y$.
\end{dfs}
Existence of suprema in the Cuntz semigroup was first observed by the second named author
in~\cite[Lemma 3.2]{P1} for C$^*$-algebras with real rank zero and stable rank one. In this
section we drop the condition of real rank zero and obtain the same result, albeit with
considerably more effort. We have been informed by George Elliott that suprema in the Cuntz
semigroup exist in full generality, a result he has proved with K. Coward and C. Ivanescu.  No
preprint was available at the time of writing, but we state for the record their result predates
ours.  It is not clear whether their result will suffice for our application, as we require a
particular description of suprema in the Cuntz semigroup.

% In order to ease notation, we let $A_a$ denote the hereditary C$^*$-algebra
% generated by a positive element $a$ in $A$, i.e.,
% $A_a=\overline{aAa}$. Recall that if $A$ is a separable
% C$^*$-algebra, then all hereditary algebras are of this form. Note
% that if $A_a\subseteq A_b$ for positive elements $a$ and $b$ in
% $A$, then we have that $a\precsim b$ (by~\cite[Proposition 2.7
% (i)]{KR}).

\begin{lms}
\label{lem:incrhereditary} Let $A$ be a unital and separable
C$^*$-algebra, and let $a_n$ be a sequence of positive elements in
$A$ such that $A_{a_1}\subseteq A_{a_2}\subseteq \cdots$. Let
$A_{\infty}=\overline{\cup_{n=1}^{\infty}A_{a_n}}$, and let
$a_{\infty}$ be a strictly positive element of $A_{\infty}$. Then
\[
\langle a_{\infty}\rangle =\sup_n \langle a_n\rangle\,.
\]
Moreover, for any trace $\tau$ in $\TA$, we have
$\dt(a_{\infty})=\sup_n\dt(a_n)$.
\end{lms}
\begin{proof}
To prove that $\langle a_{\infty}\rangle\geq \langle a_n\rangle$,
it suffices to prove that, as observed above,
$A_{\infty}=A_{a_{\infty}}$. For this, it is enough to show that
$A_{\infty}$ is hereditary. Indeed, if $a\in A$ and $c_1$, $c_2\in
A_{\infty}$, then choose sequences $x_n$ and $y_n$ in $A_{a_n}$
such that
\[
\Vert x_n-c_1\Vert\rightarrow 0\text{ and }\Vert
y_n-c_2\Vert\rightarrow 0\,.
\]
Then $x_nay_n\in A_n$, and since $c_1ac_2=\lim\limits_n x_nay_n$,
we see that $c_1ac_2$. (Recall from, e.g.~\cite[Theorem 3.2.2]{murphy}, that a
C$^*$-subalgebra $C$ of $A$ is hereditary if and only if
$c_1ac_2\in C$ whenever $a\in A$ and $c_1$, $c_2\in C$.)

Now assume that $\langle a_n\rangle\leq \langle b\rangle$ for all
$n$ in $\mathbb{N}$. Choose positive elements $x_n$ in $A_{a_n}$
such that $\Vert x_n-a_{\infty}\Vert<\delta_n $, where
$\delta_n\rightarrow 0$. It then follows by~\cite[Lemma 2.5
(ii)]{KR} that $\langle (a_{\infty}-\delta_n)_+ \rangle\leq
\langle x_n \rangle\leq \langle a_n\rangle\leq\langle b\rangle$.
Thus~\cite[Proposition 2.6]{KR} (or~\cite[Proposition
2.4]{Rfunct}) entails $\langle a_{\infty}\rangle\leq \langle
b\rangle$, as desired.

Also, since $\langle x_n\rangle\leq\langle a_n\rangle\leq\langle
a_{n+1}\rangle\leq\langle a_{\infty}\rangle$ for all $n$ and
$\lim_n x_n=a_{\infty}$, we have that, if $\tau\in\TA$,
\[
\sup_{n\to\infty}\dt(a_n)\leq\dt(a_{\infty})\leq\liminf_{n\to\infty}\dt(x_n)\leq\liminf_{n\to\infty}\dt(a_n)=\sup_{n\to\infty}\dt(a_n)
\]
\end{proof}

We shall assume in the results below that $\mathrm{sr}(A)=1$. Recall that, under this assumption, Cuntz
subequivalence is implemented by unitaries (by condition (iv) in Proposition~\ref{basics}).
% \begin{props}[R{\o}rdam, \cite{Rfunct}, Proposition 2.4]
% Let $A$ be a unital C$^*$-algebra with stable rank one, and let
% $a,b \in A_+$.  Then, $a \precsim b$ if and only if for every
% $\epsilon>0$ there exists $u \in U(A)$ such that $u^*(a-\epsilon)_+u\in
% A_b$.
% \end{props}
Note that, in this case, $a\precsim b$ implies that for each
$\epsilon>0$, there is $u$ in $U(A)$ such that
$A_{(a-\epsilon)_+}\subseteq uA_bu^*$. Indeed, if $a\in
A_{(a-\epsilon)_+}$, then find a sequence $(z_n)$ in $A$ such that
$a=\lim_n (a-\epsilon)_+z_n(a-\epsilon)_+$. Writing
$(a-\epsilon)_+=uc_{\epsilon}u^*$, with $c_{\epsilon}$ in $A_b$,
we see that $a=u(\lim_n c_{\epsilon}u^*z_n uc_{\epsilon})u^*\in
uA_b u^*$.
\begin{lms}
\label{lem:supremainA} Let $A$ be a unital and separable
C$^*$-algebra with $\mathrm{sr}(A)=1$. Let $(a_n)$ be a sequence
of elements in $A$ such that $\langle a_1\rangle\leq\langle
a_2\rangle\leq\cdots$. Then $\sup_n\langle a_n\rangle$ exists in
$\WA$ and for any $\tau$ in $\TA$, we have
$\dt(\sup_n\langle a_n\rangle)=\sup_n\dt(a_n)$.
\end{lms}
\begin{proof}
Define numbers $\epsilon_n>0$ recursively. Let $\epsilon_1=1/2$,
and choose $\epsilon_n<1/n$ such that
\[
(a_j-\epsilon_j/k)_+\precsim (a_n-\epsilon_n)_+
\]
for all $1\leq j<n$ and $1\leq k\leq n$. (This is possible
using~\cite[Proposition 2.6]{KR} and because $a_j\precsim a_n$ for
$1\leq j< n$. Notice also that $(a_n-\epsilon)_+\leq
(a_n-\delta)_+$ whenever $\delta\leq\epsilon$.)

Since $(a_1-\epsilon_1/2)_+\precsim (a_2-\epsilon_2)_+$ and
$\mathrm{sr}(A)=1$, there is a unitary $u_1$ such that
\[
A_{(a-\epsilon_1/2)_+-\epsilon_1/2)_+}\subseteq
u_1A_{(a_2-\epsilon_2)_+}u_1^*\,.
\]
But $(a-\epsilon_1/2)_+-\epsilon_1/2)_+=(a_1-\epsilon_1)_+$
(see~\cite[Lemma 2.5]{KR}), so
\[
A_{(a-\epsilon_1)_+}\subseteq u_1A_{(a_2-\epsilon_2)_+}u_1^*\,.
\]
Continue in this way, and find unitaries $u_n$ in $A$ such that
\[
A_{(a-\epsilon_1)_+}\subseteq
u_1A_{(a_2-\epsilon_2)_+}u_1^*\subseteq
\]
\[
\subseteq u_1u_2A_{(a_3-\epsilon_3)_+}u_2^*u_1^*\subseteq \cdots
\subseteq
(\prod\limits_{i=1}^{n-1}
u_i)A_{(a_n-\epsilon_n)_+}(\prod\limits_{i=1}^{n-1}u_i)^*\subseteq\cdots
\]
Use Lemma~\ref{lem:incrhereditary} to find a positive element
$a_{\infty}$ in $A$ such that
\[
\langle a_{\infty}\rangle=\sup_n \langle
(a-\epsilon_n)_+\rangle\,,
\]
and also $\dt(a_{\infty})=\sup_n\dt((a-\epsilon_n)_+)\leq\sup_n\dt
(a_n)$ for any $\tau$ in $\TA$.

We claim that $\langle a_{\infty}\rangle=\sup_n\langle a_n\rangle$
as well. From this it will readily follow that
$\dt(a_{\infty})=\sup_n\dt(a_n)$.

To see that $\langle a_n\rangle\leq \langle a_{\infty}\rangle$ for
all $n$ in $\mathbb{N}$, fix $n<m$ and recall that, by
construction,
\[
\langle (a_n-\epsilon_n/(m-1))_+\rangle\leq \langle
(a_m-\epsilon_m)_+)\rangle\leq \langle a_{\infty}\rangle\,.
\]
Hence, letting $m\rightarrow \infty$, we see that $\langle
(a_n-\epsilon)_+\rangle\leq \langle a_{\infty}\rangle$ for any
$\epsilon>0$, and so $\langle a_n\rangle\leq \langle
a_{\infty}\rangle$ for all $n$. Conversely, if $\langle
a_n\rangle\leq \langle b\rangle$ for all $n$ in $\mathbb{N}$, then
also $\langle (a_n-\epsilon_n)_+\rangle\leq \langle b\rangle$ for
all natural numbers $n$, and hence $\langle a_{\infty}\rangle\leq
\langle b\rangle$.
\end{proof}

\begin{thms}
\label{thm:cuntzsuprema} Let $A$ be a unital and separable
C$^*$-algebra with stable rank one. Then every bounded sequence
$\{\langle a_n\rangle\}$ in $\WA$ has a supremum $\langle
a_{\infty}\rangle$ and $\dt(a_{\infty})=\sup_n\dt(a_n)$ for any
$\tau$ in $\TA$.
\end{thms}
\begin{proof}
Let $\langle x_1\rangle\leq\langle x_2\rangle\leq\cdots$ be given,
and assume that $\langle x_n\rangle\leq k\langle 1_A\rangle$ for
all $n$.

Inspection of the proof of Lemma~\ref{lem:supremainA} reveals that
we may choose a sequence $\epsilon_n>0$ with the following
properties:
\begin{enumerate}
\item[(i)] $\langle (x_n-\epsilon_n)_+ \rangle \leq \langle
(x_{n+1}-\epsilon_{n+1})_+\rangle$.

\item[(ii)] If $\langle (x_n-\epsilon_n)_+\rangle\leq \langle b\rangle$
for all $n$, then $\langle x_n\rangle\leq\langle b\rangle$ for all
$n$.
\end{enumerate}
Since $\langle x_n\rangle\leq k\langle 1_A\rangle$, find $y_n$ in
$M_{\infty}(A)_+$ such that
\[
(x_n-\epsilon_n)_+=y_n(1_A\otimes 1_{M_k})y_n^*\,.
\]
Define $a_n=(1_A\otimes 1_{M_k})y_n^*y_n(1_A\otimes 1_{M_k})$,
which is an element of $M_k(A)$. Then $\langle a_n\rangle=\langle
(x_n-\epsilon_n)_+\rangle\leq \langle a_{n+1}\rangle$ for all $n$.
Since $M_k(A)$ also has stable rank one, we may use
Lemma~\ref{lem:supremainA} to conclude that $\{\langle
a_n\rangle\}$ has a supremum $\langle a_{\infty}\rangle$ with
$a_{\infty}$ in $M_k(A)$. It follows that then $\langle
a_{\infty}\rangle$ is the supremum of $\{\langle
(x_n-\epsilon_n)_+\rangle\}$ in $\WA$. Evidently, our selection of
the sequence $\epsilon_n>0$ yields that $\langle
a_{\infty}\rangle=\sup_n\langle x_n\rangle$.

The proof that $\dt(a_{\infty})=\sup_n\dt(a_n)$ is identical to
the one in Lemma~\ref{lem:supremainA}.
\end{proof}
Recall that a state $s$ on a preordered monoid $M$ with order unit
$u$ is \emph{$\sigma$-normal} if whenever $(a_n)$ is an increasing
sequence and $\sup_n a_n=a$ exists, then $s(a)=\sup_n s(a_n)$.
Denote the set of $\sigma$-normal states on $M$ by
$\mathrm{St}_{\sigma}(M,u)$.
\begin{cors}
Let $A$ be a unital, separable and exact C$^*$-algebra with stable
rank one. Then $\mathrm{LDF}(A)=\mathrm{St}_{\sigma}(\WA,\langle
1_A\rangle)$.
\end{cors}
\begin{proof}
The inclusion $\mathrm{St}_{\sigma}(\WA,\langle
1_A\rangle)\subseteq\mathrm{LDF}(A)$ always holds, as shown
in~\cite[Proposition 3.3]{P1}. The converse inclusion follows
directly from Theorem~\ref{thm:cuntzsuprema} and the fact that
every lower semicontinuous function comes from a trace
(see~\cite{bh}).
\end{proof}
\begin{cors}
Let $A$ be a unital and separable C$^*$-algebra with stable rank
one. If $x\in \WA$ is such that $x\leq \langle 1_A\rangle$, then
there is $a$ in $A$ such that $x=\langle a\rangle$.
\end{cors}
\begin{proof}
There are a natural number $n$ and an element $b$ in $M_n(A)_+$
such that $x=\langle b\rangle$. For any $m$ in $\mathbb{N}$, find
elements $x_m$ such that
\[
(b-1/m)_+=x_m1_Ax_m^*\,,
\]
so $a_m:=1_Ax_m^*x_m1_A\in A$ and $a_m\sim (b-1/m)_+$. Moreover,
the sequence $\langle a_m\rangle$ is increasing, and the proof of
Lemma~\ref{lem:supremainA} ensures that it has a supremum $a$ in
$A$. Clearly,
\[
\langle a\rangle=\sup_m\langle a_m\rangle=\sup_m\langle
(b-1/m)_+\rangle=\langle b\rangle\,.
\]
\end{proof}

\begin{cors}\label{cor:whensupproj}
Let $A$ be a unital and separable C$^*$-algebra with stable rank
one. If $\langle a_n\rangle$ is a bounded and increasing sequence
of elements in $W(A)$ with supremum $\langle a\rangle$. Then
$\langle a\rangle=\langle p\rangle$ for a projection $p$, if and
only if, there exists $n_0$ such that $\langle a_n\rangle=\langle
p\rangle$ whenever $n\geq n_0$.
\end{cors}
\begin{proof}
Suppose that $\langle a\rangle=\sup_n\langle a_n\rangle=\langle
p\rangle$ for a projection $p$. We may assume that all the
elements $a$, $a_n$ and $p$ belong to $A$. For any $n$, we have
that $a_n\precsim p$. On the other hand, the proof of
Lemma~\ref{lem:supremainA} shows that $p=\lim_n b_n$, for some
elements $b_n\precsim (a_n-\epsilon_n)_+$ (where $\epsilon_n>0$ is
a sequence converging to zero). From this it follows that for
sufficiently large $n$, $p\precsim b_n\precsim
(a_n-\epsilon_n)_+\precsim a_n$. Thus $p\sim a_n$ if $n$ is large
enough, as desired.
\end{proof}

%%%%%%%%%%%%%%%%%%%%%%%%%%%%%%%%%%%%%%%%%%%
\section{Surjectivity of $\iota\colon \WA_+ \to \lap$}
\label{sec:surjectivity}

In this section we will prove the surjectivity of $\iota$ for algebras satisfying
the hypotheses of Theorems A and B.  This will complete the proof of Theorem A.

Our first proposition follows from Corollary \ref{cor:realizingaffinefunctions}.  In a break with convention,
we let $\mathrm{CAff}(\bullet)$ denote {\it continuous} affine functions for the remainder of the
paper --- this is necessary as we deal also with not-necessarily-continuous affine functions.

\begin{props}\label{prop:denseunion}  Assume $A_1 \subset A_2 \subset \cdots \subset A$
are unital subalgebras with dense union.  If $A$ is simple and $f \in \baff$ is continuous
and strictly positive, then for every $\varepsilon > 0$ there exists $n \in \mathbb{N}$
and $0 \leq a \in A_n$ such that $|f(\tau) - \tau(a)| < \varepsilon$, for all $\tau \in \TA$.
(Using self-adjoint $a$, this holds without simplicity.)

Consequently, there exists a continuous function $g \in \mathrm{Aff}(\mathrm{T}(A_n))$ ---
namely, $\hat{a}$ --- whose image under the canonical map $\mathrm{Aff}(\mathrm{T}(A_n))
\to \baff$ is within $\varepsilon$ of $f$.
\end{props}

\begin{lms}\label{approx}
Let $A = p(\mathrm{C}(X) \otimes \mathcal{K})p$ be a homogeneous C$^*$-algebra with
$X$ a compact metric space and $\mathrm{rank}(p)=n$.  Let there be given $g \in \mathrm{CAff}(\TA)$ satisfying
$0 \leq g \leq 1$.  Then, there exists $a \in \mathrm{M}_{\infty}(A)_+$
such that $f := \iota(a)$ satisfies
\[
0 \leq g(\tau) - f(\tau) \leq 1/n, \ \forall \tau \in \TA\,.
\]
\end{lms}

\begin{proof}
For each $0 \leq i \leq n-1$ define an open set
\[
A_i := \{ x \in X \ | \ g(x) > i/n \}.
\]
Notice that $A_i \subseteq A_j$ whenever $j \leq i$.  Since $X$ is metric, we
can find, for each $i$, a continuous function $f_i:X \to [0,1]$ such that
$f_i(x) \neq 0$ if and only if $x \in A_i$.  Put
\[
B_i := \{ x \in X \ | \ (i+1)/n \geq g(x) > i/n \} = A_i \backslash (\cup_{j > i} A_j).
\]
and
\[
a := \bigoplus_{i=1}^{n-1} f_i \cdot q \in \mathrm{M}_{\infty}(A)_+,
\]
where $q$ is a fixed rank one projection in some $\mathrm{M}_n(\mathrm{C}(X)) \subseteq \mathrm{M}_{\infty}(A)$.

The tracial simplex of $A$ is a Bauer simplex, so the lower semicontinuous
affine functions on $\TA$ are in bijective correspondence with the
lower semicontinuous functions on the extreme boundary $\partial_e \TA \cong X$ via
restriction.
For each $x \in X$, the value of $f(x) := \iota(a)(x)$ is the normalised rank of $a$ at $x$.
In other words,
\[
\iota(a)(x) := \frac{| \{ j \geq 1 \ | \ x \in A_j \} |}{n}.
\]
If $x \in (X \backslash A_0) \cup B_0$, then $f(x) = 0$, and $0 \leq (g-f)(x) \leq 1/n$
for all such $x$.  If $j \geq 1$ and $x \in B_j$, then $f(x) = j/n$ and
$j/n < g(x) \leq (j+1)/n$, and $0 \leq (g-f)(x) \leq 1/n$ for all such $x$.
Since $f$ is lower semicontinuous, so is $f-g$.  A lower semicontinuous affine
function on a Bauer simplex achieves its minimum on the extreme boundary, and this
minimum is at least $-1/n$ by construction.  Thus, $f-g \geq -1/n$.  By affineness,
$f-g \leq 0$ on every finite convex combination of extreme traces.  Every point
$\tau \in \TA$ is the weak-$*$ limit of a sequence of such combinations,
so the lower semicontinuity of $f-g$ yields $f-g \leq 0$ on $\TA$.
\end{proof}

Let $A$ be a unital C$^*$-algebra.
It is well known that $A \mapsto \mathrm{CAff}(\TA)$ is a covariant functor
into the category of complete order-unit spaces.   If $B$ is
a unital C$^*$-algebra and $\psi:A \to B$ is a $*$-homomorphism, then
let
\[
\psi^{\sharp}:\mathrm{T}(B) \to \TA
\]
denote the map induced on traces.  The induced map
\[
\psi^{\bullet}:\mathrm{CAff}(\TA) \to \mathrm{CAff}(\mathrm{T}(B))
\]
is then given by
\[
\psi^{\bullet}(f)(\gamma) = f(\psi^{\sharp}(\gamma)).
\]
Let $a \in A$ be positive, with image $\iota(a) \in \mathrm{LAff}_b(\TA)^+$.
Then, $\iota(\psi(a)) = \psi^{\bullet}(\iota(a))$.  Indeed, for $\gamma \in \mathrm{T}(B)$
we have
\begin{eqnarray*}
\iota(\psi(a))(\gamma) & = & \lim_{n \to \infty} \gamma(\psi(a)^{1/n}) \\
& = & \lim_{n \to \infty} \gamma(\psi(a^{1/n})) \\
& = & \lim_{n \to \infty} \psi^{\sharp}(\gamma)(a^{1/n}) \\
& = & \iota(a)(\psi^{\sharp}(\gamma)) \\
& = & \psi^{\bullet}(\iota(a))(\gamma)
\end{eqnarray*}

%Let $(A_i,\phi_i)_{i \in \mathbb{N}}$ be a direct system of unital C$^*$-algebras
%with limit $A$, and let $\phi_{i\infty}:A_i \to A$ denote the canonical map.
%The construction of the inductive limit order-unit space $\mathrm{CAff}(\TA)$
%shows that the set
%\begin{equation}\label{posdense}
%\bigcup_{i =1}^{\infty} \phi_{i\infty}^{\bullet}(\mathrm{CAff}(\mathrm{T}(A_i))^+)
%\end{equation}
%is uniformly dense in $\mathrm{CAff}(\TA)^+$, where the superscript ``+'' indicates
%non-negative functions.

\begin{thms}\label{thm:ahsurjective}
Let $A$ be a simple, unital, separable, and infinite-dimensional AH algebra of stable rank
one.  If $A$ has strict comparison of positive elements, then the map
\[
\iota:W(A)_+ \to \mathrm{LAff}_b(\TA)^{++}
\]
is surjective.
\end{thms}

\begin{proof}

By Theorem \ref{thm:cuntzsuprema} and Corollary \ref{cor:whensupproj}
 it will suffice to find, for any $f \in \mathrm{LAff}_b(\TA)^{++}$,
a sequence $(a_i)_{i=1}^{\infty}$ in $A_{+}$ such that $a_i \precsim a_{i+1}$, $\langle a_i \rangle
\neq \langle a_{i+1} \rangle$, and
\[
\lim_{i \to \infty} d_{\tau}(a_i) = f(\tau).
\]

First, use the lower semicontinuity of $f$ to find a sequence $(f_i)_{i=1}^{\infty}$
in $\mathrm{CAff}(\TA)^{++}$ satisfying
\begin{enumerate}
\item[(i)] $f_i(\tau) < f_{i+1}(\tau)$ for every $i \in \mathbb{N}$
and $\tau \in \TA$, and
\item[(ii)] $\lim_{i \to \infty} f_i(\tau) = f(\tau)$ for every $\tau \in \TA$.
\end{enumerate}
Since the difference $f_{i}-f_{i-1}$ is continuous and strictly positive on the
compact space $\TA$, it achieves a minimum, say $\epsilon_i>0$.

Let $A = \lim_{i \to \infty}(A_i,\phi_i)$ be an AH decomposition for $A$, i.e.,
\[
A_i = \bigoplus_{j=1}^{n_i} p_{i,j}(\mathrm{C}(X_{i,j}) \otimes \mathcal{K})p_{i,j}
\]
for compact connected metric spaces $X_{i,j}$ and projections $p_{i,j} \in \mathrm{C}(X_{i,j}) \otimes \mathcal{K}$.
Put $A_{i,j} := p_{i,j}(\mathrm{C}(X_{i,j}) \otimes \mathcal{K})p_{i,j}$.
By Proposition \ref{prop:denseunion} we may assume, modulo compression of our inductive
system, that $f_i \in  \phi_{i\infty}^{\bullet}(\mathrm{CAff}(\mathrm{T}(A_i))^+)$ for
each $i \in \mathbb{N}$.  Let $\tilde{f}_i$ be a pre-image of $f_i$ in $\mathrm{CAff}(\mathrm{T}(A_i))^{+}$.
By compressing our inductive sequence again if necessary we may, by the simplicity and non-finite-dimensionality
of $A$, assume that
\[
\frac{1}{\mathrm{min}_j \ \mathrm{rank}(p_{i,j})} \ll \epsilon_i.
\]
Use Lemma \ref{approx} to find, for each $1 \leq j \leq n_i$, an $a_{i,j} \in
\mathrm{M}_{\infty}(A_{i,j})_+$ such that
\[
0 \leq \tilde{f}_i|_{A_{i,j}} - \iota(a_{i,j}) \leq \epsilon_i/2.
\]
Put $\tilde{a}_i := \sum_{j=1}^{n_i} a_{i,j}$.  Then,
\[
0 \leq \tilde{f}_i - \iota(\tilde{a}_{i}) \leq \epsilon_i/2.
\]
The inequalities above are preserved under $\phi_{i \infty}^{\bullet}$, so that
with $a_i := \phi_{i \infty}(\tilde{a}_i)$ we have
\[
0 \leq f_i - \iota(a_i) \leq \epsilon_i/2.
\]
One easily checks that $\lim_{i \to \infty} d_{\tau}(a_i) = f(\tau)$ for each
$\tau \in \TA$.  Moreoever, we have $\iota(\langle a_i\rangle) < \iota(\langle a_{i+1} \rangle)$,
whence $\langle a_i \rangle \neq \langle a_{i+1} \rangle$ and $a_i \precsim a_{i+1}$.
\end{proof}

Now we consider the $\Z$-stable case.

\begin{lms} Let $X$ be a compact metric space and $f \in \mathrm{Aff}(\mathrm{T}(C(X)\otimes Z))$ be a nonnegative lower semicontinous function.  Then, there exists an element $\langle a \rangle  \in \mathrm{W}(C(X)\otimes Z)$ such that $\| \iota(\langle a\rangle) - f \| < \epsilon$.
\end{lms}

\begin{proof}  Since the tracial simplex of $C(X)\otimes Z$ is affinely homeomorphic to that of $C(X)$, we are again in the situation of a Bauer simplex.  We first handle the case that $f = \chi_{\mathcal{O}}$, where $\mathcal{O} \subset X$ is an open set.  As before, just define $a \in C(X)$ to be any function which is positive precisely on $\mathcal{O}$ and one has $\iota(\langle a \rangle) = \chi_{\mathcal{O}}$.

We can even hit multiples of such characteristic functions. Indeed, if $0 < t < 1$ we can find an element $z_t \in \Z$ such that $\iota(a\otimes z_t)$ equals $t$ times $\chi_{\mathcal{O}}$ (cf.\ \cite[Proposition 3.2]{pt}).  This, however, completes
the proof since linear combinations of such characteristic functions are uniformly dense in the
lower semicontinuous functions.
\end{proof}

\begin{thms} Let $A$ be any simple, unital, and exact C$^*$-algebra which is finite and $\Z$-stable.
Then,
\[
\iota\colon\WA_+ \to \mathrm{LAff}_b(\TA)^{++}
\]
is surjective.
\end{thms}

\begin{proof}  It suffices to show that if $f \in \mathrm{LAff}_b(\TA)^{++}$ is continuous
then we can approximate it arbitrarily well by elements in $\iota(\WA_+)$.

By Corollary \ref{cor:realizingaffinefunctions}, we can find $0 \leq a \in A$ such that $f = \hat{a}$.
Let $\psi:\Z \otimes \Z \to \Z$ be any $*$-isomorphism, and define
\[
\phi:A \otimes \Z \otimes \Z \to A \otimes \Z \otimes \Z
\]
by
\[
\phi(a \otimes z_1 \otimes z_2) = a \otimes \psi(z_1 \otimes z_2) \otimes 1_{\Z}.
\]
By \cite[Corollary 1.12]{TW1}, $\phi:A \to A$ is approximately inner, whence $\widehat{\phi(a)} = \hat{a}$.
We will thus assume below that upon identifying $A$ with $A \otimes \Z$, we have $a \in A \otimes 1_{\Z}$.

Let $B = C^*(a) \otimes \Z$ and now regard $\hat{a}$ as a continuous affine function on the tracial space
of $B$.  By the previous lemma we can approximate $\hat{a} \in \mathrm{Aff}_b(\mathrm{T}(B))$ by the image
of $\mathrm{W}(B)$. By functoriality, it follows that $f$ is approximated by $\iota(\WA_+)$.
\end{proof}

\begin{rems} {\rm It is proved in \cite{tomscomp} that a simple, unital, and infinite-dimensional
AH algebra of slow dimension growth has strict comparison;  such algebras also have stable rank one
by the main results of \cite{bdr}. }
\end{rems}

\begin{cors}\label{isomorphism} Let $A$ be a simple, unital, and finite
 C$^*$-algebra which is either exact and $\Z$-stable or an
infinite-dimensional AH algebra of slow dimension growth.  Then,
\[
\phi:\WA \to \widetilde{\mathrm{W}}(A)
\]
is an order isomorphism
\end{cors}

\begin{proof} Knowing the surjectivity of $\iota$ for these two classes of algebras, the result follows
from Theorem \ref{embed}.
\end{proof}

We conjecture that Corollary \ref{isomorphism} holds for simple, separable, unital ASH algebras
with strict slow dimension growth, and so, by deep results of Q. Lin and N. C. Phillips, for
a large class of C$^*$-dynamical systems.

%%%%%%%%%%%%%%%%%%%%%%%%%%%%%%%%%%%%%%%%%%%%%
\section{The Blackadar-Handelman Conjectures}
\label{sec:BHconjectures}

In this section we prove Theorem B of the introduction.  In fact, we prove the Blackadar-Handelman
conjectures in somewhat greater generality.
Throughout this section $\phi:\WA \to \widetilde{\mathrm{W}}(A)$ is the map defined in section
\ref{sec:surjectivity}.

\begin{lms} \label{lem:density} Let $\mathcal{S} \subset \baff$ be any sub-semigroup
containing the constant function $1$, endowed with the pointwise
(pre)order.  If $\varphi\colon \mathcal{S} \to \mathbb{R}$ is any
state then there exists a net of traces $\{\tau_{\lambda}\}_{\lambda
\in \Lambda} \subset \TA$ such that $$\varphi(s) = \lim_{\lambda \to
\infty} s(\tau_{\lambda}),$$ for all $s \in \mathcal{S}$.
\end{lms}

\begin{proof} Thanks to \cite[Corollary 2.7]{br}, we may extend the
state $\varphi$ to a state on all of $\baff$; i.e., we may assume
$\mathcal{S} = \baff$.

However, every state on $\baff$ is actually a bounded linear functional (cf.\ \cite[Lemma 6.7]{G}).
That is, $\varphi \in (\baff)^* = Z_{sa}^*$, by Lemma \ref{lem:Zsa-spanTA}. Moreover,
$\varphi$ defines a positive linear functional on $Z_{sa}^*$, since $\varphi(0) = 0$ and $\varphi$
is order preserving.  Since the normal states on $Z$ are weak-$*$ dense in the set of all states,
it follows that $\varphi \in Z_{sa}^*$ can be approximated by a net $\{\tau_{\lambda}\}_{\lambda
\in \Lambda} \subset \TA$.
\end{proof}

The following lemma is well known.

\begin{lms}\label{infseq}
Every infinite-dimensional $C^*$-algebra contains a positive element with
infinite spectrum.
\end{lms}

\begin{cors}\label{purepos}
Let $A$ be a simple, unital, and infinite-dimensional $C^*$-algebra.
Then, $A$ contains a purely positive element.
\end{cors}

\begin{proof}
By the previous lemma, there is a positive element $a \in A$ with infinite
spectrum.  Choose an accumulation point $x \in \sigma(a)$.  Let $f$ be a continuous
function on $\sigma(a)$ such that $f(t)$ is nonzero if and only if $t \neq x$.
Then, $f(a)$ is positive and has zero as an accumulation point of its spectrum.
$f(a)$ is thus purely positive by \cite[Proposition 2.1]{pt}.
\end{proof}

\begin{thms}
\label{thm:bh2} Let $A$ be a simple, unital, exact, and stably finite C$^*$-algebra
for which
\[
\phi:\WA \to \widetilde{\mathrm{W}}(A)
\]
is an order-embedding.  Then,
$\mathrm{LDF}(A)$ is dense in $\mathrm{DF}(A)$.
\end{thms}

\begin{proof}
We may assume that $A$ is infinite-dimensional, whence $\WA_+$ is non-empty by Corollary \ref{purepos}. Thus,
$\mathrm{K}_0^*(A)$ is order-isomorphic to $G(\WA_+)$
(see~\cite[Lemma 5.5]{pt}). Let $\gamma\colon \WA_+\to
G(\WA_+)$ denote the natural Grothendieck map.

If we pick any $c$ in $\WA_+$, then we can define an
order-isomorphism $\alpha$ by
\[
\alpha ([p])=\gamma(\langle p\rangle +c)-\gamma (c)
\]
if $p$ is a projection, and
\[
\alpha([a])=\gamma(\langle a\rangle)
\]
if $\langle a\rangle\in \WA_+$.
We thus have that, by composition, $\mathrm{K}_0^*(A)$ is order-isomorphic
to a subgroup $\mathcal{S}$ of $\{f-g\mid f,g\in\laff_b
(\TA)^{++}\}$ via
\[
[a]-[b]\mapsto \widehat{a}-\widehat{b}\,,
\]
where $\widehat{a}(\tau)=d_\tau(a)$ (for any $\tau$ in
$\TA$). Note that under this order-isomorphism, $[1]$
is mapped to $(1\oplus
c')^{\widehat{}}-\widehat{c'}=1+\widehat{c'}-\widehat{c'}=1$,
where $c'$ is any purely positive element such that $\langle
c'\rangle=c$.

Next, if $d\in \mathrm{DF}(A)$, then by the isomorphism we may think of $d$
as a normalized state on the image $\mathcal{S}$, which is a
subsemigroup of $\baff$ containing the constant function $1$. By
Lemma~\ref{lem:density}, there is a net of
traces $\{\tau_\lambda\}$ in $\TA$ such that
$d(s)=\lim\limits_{\lambda}s(\tau_\lambda)$ for any $s$ in
$\mathcal{S}$. In particular, for $a$ in $A$:
\[
d([a])=\lim\limits_{\lambda}(\widehat{a}(\tau_\lambda))=\lim\limits_{\lambda}d_{\tau_\lambda}(a)\,,
\]
and since $a\mapsto d_{\tau_\lambda}(a)$ is in $\mathrm{LDF}(A)$, the proof
is complete.
\end{proof}

\begin{rems} {\rm
The order-embedding hypothesis above is satisfied whenever $A$ has strict comparison. For
example, it suffices to know $A$ is $\mathcal{Z}$-stable or an AH algebra of slow dimension growth,
though this is overkill as it implies $\phi$ is an order-isomorphism. }
\end{rems}

\begin{dfs}
Let $(M,\leq)$ be a preordered monoid. We say that $M$ satisfies
the \emph{Riesz Interpolation Property} if whenever $x_1$, $x_2$,
$y_1$, $y_2\in M$ satisfy $x_i\leq y_j$ for all $i$ and $j$, then
there is $z$ in $M$ such that $x_i\leq z\leq y_j$.
\end{dfs}

\begin{lms}
\label{lem:interpolation} Let $K$ be a metrizable compact convex
set. Then $\laff_b(K)^{++}$, equipped with the pointwise ordering,
is an interpolation monoid.
\end{lms}
\begin{proof}
Let there be given functions $f_1$, $f_2$, $g_1$, $g_2$ in
$\laff_b(K)^{++}$ such that $f_i\leq g_j$ for $i$, $j=1$, $2$.

Since $K$ is metrizable, we may write $f_i=\sup_n f_{i,n}$, where
$f_{i,n}\in \mathrm{CAff}(K)^{++}$ and $f_{i,n}\leq f_{i,n+1}$ for $i=1$,
$2$ and all $n$. There is $h_1$ in $\mathrm{CAff}(K)^{++}$ such that
$f_{i,1}\leq h_1\leq g_j$, by, e.g.~\cite[]{G}.

Next, since $f_{i,2},h_1\leq g_j$ ($i$, $j=1$, $2$), there is
$h_2$ in $\mathrm{CAff}(K)^{++}$ such that
\[
f_{i,2},h_1\leq h_2\leq g_j\,.
\]

Continue in this way to find an increasing sequence $h_n$ in
$\mathrm{CAff}(K)^{++}$ such that $f_{i,n}\leq h_n\leq g_j$ for $i$, $j=1$,
$2$ and all $n$. Put $h=\sup_n h_n$, which is an element of $\laff
(K)^{++}$ (as it is a supremum of continuous and affine
functions). Then, by construction $f_i\leq h\leq g_j$ for all $i$,
$j$.
\end{proof}
\begin{thms}
\label{thm:bhconj1} Let $A$ be a simple, unital, exact, and stably finite
C$^*$-algebra. If
\[
\phi:\WA \to \widetilde{\mathrm{W}}(A)
\]
is an
order isomorphism, then $\mathrm{DF}(A)$ is a Choquet simplex.
\end{thms}
\begin{proof}
We may assume that $A$ is infinite-dimensional --- the finite-dimensional
case follows from the fact that $\VA \cong \WA$ (\cite{tomscomp}).

Since $A$ is infinite dimensional, the semigroup $\WA_+$ is
non-empty by Corollary \ref{purepos}. Thus, we may
use~\cite[Lemma 5.2]{pt}, which ensures that the partially ordered
group $K_0^*(A)$ is order-isomorphic to $G(\WA_+)$ (with its
natural ordering induced by the partial order in $\WA_+$). Since,
as just mentioned, $\WA_+\cong\laff_b(\TA)^{++}$,
Lemma~\ref{lem:interpolation} applies to conclude that $\WA_+$ is
an interpolation monoid. But then we can use~\cite[Lemma 4.2]{P1},
to see that $G(\WA_+)$ is an interpolation group.

Therefore, $(K_0^*(A), K_0^*(A)^{++})$ is an interpolation group
and thus $\mathrm{DF}(A)$, being the state space of $K_0^*(A)$, is
a Choquet simplex, by e.g.~\cite[Theorem 10.17]{G}.
\end{proof}

Combining Theorems \ref{thm:bh2} and \ref{thm:bhconj1} with Corollary \ref{isomorphism} now
yields Theorem B.

\bibliographystyle{amsplain}

\providecommand{\bysame}{\leavevmode\hbox
to3em{\hrulefill}\thinspace}

\end{document}